\documentclass{article}
\usepackage{amssymb}
\usepackage{amsfonts}
\usepackage{amsmath}

\newtheorem{theorem}{Theorem}

\newtheorem{corollary}[theorem]{Corollary}

\newtheorem{definition}[theorem]{Definition}

\newtheorem{proposition}[theorem]{Proposition}

\input{tcilatex}
\begin{document}

\title{Horizontal cohomology\\
of a local Lie group}
\author{Erc\"{u}ment Orta\c{c}gil}
\maketitle

\begin{abstract}
We define and study the invariant linear and nonlinear horizontal double
complexes of a local Lie group.
\end{abstract}

\section{Introduction}

Let $M$ be a smooth manifold with $\dim M=n$ and $\mathcal{F}\rightarrow M$
be fiber bundle with $k$'th jet extension $J_{k}\mathcal{F}\rightarrow M.$ A
fibered submanifold $\mathcal{E\subset }$ $J_{k}\mathcal{F}\rightarrow M$
defines a $k$'th order $PDE$ on $M.$ The horizontal cohomology of $\mathcal{R%
}$ is defined in [17] as a part of the general formalism of Vinagradov
spectral sequence and studied further in [18], [14], [15], [16] and by
various authors (the references in [14] contain an extensive list on
horizontal cohomology). In the presence of a group structure as in this
work, it is possible to define also the invariant horizontal complex as the
\textquotedblleft edge sequence\textquotedblright\ of the invariant
variational bicomplex (see [2], [9]).

The study of local Lie groups is initiated in [11], [1], [8], [13]. As shown
in [11], the theory of local Lie groups is not a simple consequence of the
global theory but has its own set of interesting and delicate geometric
structures. Slightly modifying the definition of a local Lie group in [11],
we showed in [1] that a Lie group can be defined as a globalizable local Lie
group, hence reinstating the paradigm of local to global to its historical
record. As we indicated in [12], local Lie groups are particular
pre-homogeneous geometric structures with vanishing curvatures (see also [3]
for a similar approach to geometric structures based on Cartan algebroids).\

A local Lie group is defined by a first order nonlinear $PDE$ on $%
J_{1}(M\times M)$ and its Lie algebra as a first order linear $PDE$ on $%
J_{1}(T)$ where $T\rightarrow M$ is the tangent bundle. The elementary
nature of these $PDE$'s allows one to make a concrete study of their
invariant horizontal cohomologies which is the purpose of this work
(Sections 1, 2). We show that the group structure enables one to realize
these horizontal complexes as the second rows of two double complexes
(Section 4). The linearization map determines a homomorphism from the
nonlinear double complex to the linear one (Sections 3, 4). In the nonlinear
case, it turns out that the first row of this invariant double complex
computes the Lie algebra cohomology and the first coloumn computes the Lie
group cohomology in analogy with the Van Est spectral sequence. In
particular, the nonlinear double complex assigns an infine number of
seemingly new cohomology groups to a Lie group and uncovers some direct
links between Lie group cohomology as initiated in [6], [7] and later
generalized to groupoids in [19], [5] and the horizontal cohomology
mentioned in the above works.

\section{Local Lie groups}

In this section we shortly recall the theory of local Lie groups. We refer
to [1] for more details for some points (see also Section 2 of [13]).

Let $(M,\widetilde{\varepsilon })$ be a manifold with a splitting $%
\widetilde{\varepsilon }$ of $J_{1}(M\times M)\rightarrow M\times M$. So $%
\widetilde{\varepsilon }$ assigns to any ordered pair $(p,q)$ a $1$-arrow
from $p$ to $q$ and this assignment preserves the composition and inversion
of arrows. Such a splitting exists if and only if $M$ is parallelizable. We
define the components

\begin{equation}
\widetilde{\Gamma }_{kj}^{i}(x)\overset{def}{=}\left[ \frac{\partial
\widetilde{\varepsilon }_{j}^{i}(x,y)}{\partial y^{k}}\right] _{y=x}=-\left[
\frac{\partial \widetilde{\varepsilon }_{j}^{i}(x,y)}{\partial x^{k}}\right]
_{x=y}\overset{def}{=}\widehat{\Gamma }_{jk}^{i}(x)
\end{equation}

For a vector field $\xi =(\xi ^{i}),$ we define $\widetilde{\nabla }_{j}\xi
^{i}\overset{def}{=}\frac{\partial \xi ^{i}}{\partial x^{j}}-\widetilde{%
\Gamma }_{ja}^{i}\xi ^{a}$ and $\widehat{\nabla }_{j}\xi ^{i}\overset{def}{=}%
\frac{\partial \xi ^{i}}{\partial x^{j}}-\widehat{\Gamma }_{ja}^{i}\xi ^{a}.$
The actions of the covariant differentiation operators $\widetilde{\nabla },$
$\widehat{\nabla }$ extend naturally from vector fields to all tensor
fields. A tensor field $t$ is called $\widetilde{\varepsilon }$-invariant if
$\widetilde{\varepsilon }(p,q)_{\ast }t(p)=t(q),$ $\widetilde{\nabla }$%
-invariant if $\widetilde{\nabla }t=0$ and $\widehat{\nabla }$-invariant if $%
\widehat{\nabla }t=0.$ It is easy to show that $t$ is $\widetilde{%
\varepsilon }$-invariant if and only if it is $\widetilde{\nabla }$%
-invariant (Proposition 1 in [13], see also Proposition 5 below). So $%
\widehat{\nabla }$-invariance is defined without the object $\widehat{%
\varepsilon }$ whose definition needs a further assumption (see (5) below).
Since the linear $PDE$ $\widetilde{\nabla }\xi =0$ admits $\widetilde{%
\varepsilon }$-invariant vector fields as solutions, its integrability
condition $\widetilde{\mathfrak{R}}=0$ is satisfied. Let $\widehat{\mathfrak{%
R}}=0$ denote the integrability condition of $\widehat{\nabla }\xi =0.$ We
define the torsion tensors $\widetilde{T}$, $\widehat{T}$ by $\widetilde{T}%
_{jk}^{i}\overset{def}{=}\widetilde{\Gamma }_{jk}^{i}-\widetilde{\Gamma }%
_{kj}^{i}$ and $\widehat{T}_{jk}^{i}\overset{def}{=}\widehat{\Gamma }%
_{jk}^{i}-\widehat{\Gamma }_{kj}^{i}.$ Clearly $\widetilde{T}=-$ $\widehat{T}
$. Straightforward computations using the definitions prove the fundamental
formulas

\begin{eqnarray}
\widetilde{\nabla }_{i}\widetilde{T}_{kl}^{j} &=&\widehat{\mathfrak{R}}%
_{kl,i}^{j} \\
\widehat{\nabla }\widehat{T} &=&0\text{ \ \ \ }if\text{ \ }\widehat{%
\mathfrak{R}}=0  \notag
\end{eqnarray}

The splitting $\widetilde{\varepsilon }$ determines the nonlinear $PDE$

\begin{equation}
\frac{\partial f^{i}(x)}{\partial x^{j}}=\widetilde{\varepsilon }%
_{j}^{i}(x,f(x))
\end{equation}%
with the integrability condition $\ \widetilde{\mathcal{R}}=0.$

\begin{definition}
$(M,\widetilde{\varepsilon })$ is a local Lie group if $\widetilde{\mathcal{R%
}}=0.$
\end{definition}

In this case the local solutions of (3) are uniquely determined by their
initial conditions $f(p)=q$ and they form a simply transitive pseudogroup on
$M$ denoted by $\widetilde{\mathcal{G}}.$ If all $f\in \widetilde{\mathcal{G}%
}$ extend (necessarily uniquely) to global diffeomorphisms of $M$, then $(M,%
\widetilde{\varepsilon })$ is called globalizable. In this case $\widetilde{%
\mathcal{G}}$ becomes a global transformation group of $M$ which acts simply
transitively. It is a fundamental fact that $\widetilde{\mathcal{R}}%
=0\Leftrightarrow \widehat{\mathfrak{R}}=0,$ the implication $\Leftarrow $
being the Lie's third fundamental theorem. For a local Lie group $(M,%
\widetilde{\varepsilon })$, the solutions $\widetilde{\Theta }$ of $%
\widetilde{\nabla }\xi =0$ becomes a Lie algebra of vector fields on $M$
which can be localized at any point $p\in M.$ As a crucial fact, it is
\textit{not }$\widetilde{\Theta }$ that integrates to $\widetilde{\mathcal{G}%
}$ but $\widehat{\Theta }$ to be defined below.

For a local Lie group $(M,\widetilde{\varepsilon })$, let $g(a,b,z)$ denote
the unique local solution of (3) in the variable $z$ satisfying the initial
condition $a\rightarrow b.$ We fix some $p,q\in $ $(U,x^{i})$ and define

\begin{equation}
\widehat{\varepsilon }_{j}^{i}(p,q)\overset{def}{=}\left[ \frac{\partial
g^{i}(p,x,q)}{\partial x^{j}}\right] _{x=p}
\end{equation}

Choosing $f\in \widetilde{\mathcal{G}}$ with $f(p)=q$ and replacing $p,x$
with $p^{\prime },x^{\prime }$ close to $p,x,$ (4) shows that the local
diffeomorphism $h:x\rightarrow g(p,x,q)$ satisfies $h(p)=q$ and is the
unique local solution of

\begin{equation}
\frac{\partial h^{i}(x)}{\partial x^{j}}=\widehat{\varepsilon }%
_{j}^{i}(x,h(x))
\end{equation}%
satisfying the initial condition $h(p)=q$. In particular, the integrability
condition $\widehat{\mathcal{R}}=0$ of (5) is satisfied. Note that $\widehat{%
\varepsilon }(p,q)$ is defined for sufficiently close $p,q$ unless $(M,%
\widetilde{\varepsilon })$ is globalizable. We now check that $\widehat{%
\varepsilon }$ is a (local) splitting. In analogy with (3) we also check $%
\left[ \frac{\partial \widehat{\varepsilon }_{j}^{i}(x,y)}{\partial y^{k}}%
\right] _{y=x}=\widehat{\Gamma }_{kj}^{i}(x).$ The local solutions of (5)
define the locally transitive pseudogroup $\widehat{\mathcal{G}}.$ If $(M,%
\widetilde{\varepsilon })$ is globalizable, that is, if the local
transformations of $\widetilde{\mathcal{G}}$ globalize, then so do the local
transformations of $\widehat{\mathcal{G}}$. In this case we have the map $%
\Psi :\widetilde{\mathcal{G}}\rightarrow \widehat{\mathcal{G}}$ defined as
follows: let $f\in \widetilde{\mathcal{G}}$ and fix some $p\in M.$ Then $%
\Psi (f)$ is the unique transformation of $\widehat{\mathcal{G}}$ whose $1$%
-arrow from $p$ to $q$ is $\widehat{\varepsilon }(p,q).$ This definition
does not depend on $p$ and $\Psi (f\circ g)=\Psi (f)\circ \Psi (g),$ $\Psi
(f^{-1})=\Psi (f)^{-1}.$ As expected, a tensor field $t$ is $\widehat{%
\varepsilon }$-invariant if and only if it is $\widehat{\nabla }$-invariant.
Now $\widetilde{\Theta }$ integrates to $\widehat{\mathcal{G}}$ and the Lie
algebra $\widehat{\Theta }$ of solutions of $\widehat{\nabla }\xi =0$
integrates to $\widetilde{\mathcal{G}},$ that is, $\widetilde{\Theta },$ $%
\widehat{\Theta }$ are the \textquotedblleft Lie algebras\textquotedblright\
of the transformation groups $\widehat{\mathcal{G}},$ $\widetilde{\mathcal{G}%
}$ respectively. This corresponds to the well known fact that on a Lie group
left (right) invariant vector fields integrate to right (left) translations.
However, observe that there is no canonical choice of left and right for a
local Lie group even if it is globalizable. It is for this reason that we
avoid the notation $\mathcal{G}_{L}$ (or $\mathcal{G}_{R})$ for $\widetilde{%
\mathcal{G}}$. However, observe that the roles of $\ \widetilde{}$ and \ $%
\widehat{}$ \ are not symmetric unless $\widetilde{\mathcal{R}}=0$ and $(M,%
\widetilde{\varepsilon })$ is globalizable. Some contemplation reveals that
the static concepts of left/right on a Lie group are replaced with the
dynamic concept of \textquotedblleft time\textquotedblright\ in a local Lie
group. Now the isomorphism $\Psi :\widetilde{\mathcal{G}}\rightarrow
\widehat{\mathcal{G}}$ determines the isomorphism $d\Psi :\widehat{\Theta }%
\rightarrow \widetilde{\Theta }$ as follows: let $\xi \in \widehat{\Theta }$
and fix some $p\in M$ . We define $d\Psi (\xi )$ as the unique $\eta \in
\widetilde{\Theta }$ satisfying $\eta (p)=\xi (p).$ This definition is again
independent of $p.$

\section{The linear horizontal complex}

Let $(M,\widetilde{\varepsilon })$ be a local Lie group and $\pi
:T\rightarrow M$ be the tangent bundle

\begin{definition}
A horizontal $k$-form on $\pi :T\rightarrow M$ is a function which assigns
to any $\xi \in T$ an ordinary $k$-form at $\pi (\xi )\in M.$
\end{definition}

For the moment, we do not assume that a horizontal form is linear on the
fibers $\pi ^{-1}(x),$ $x\in M.$ We denote the vector space of the
horizontal $k$-forms by $\Lambda _{hor}^{k}(T)$. Choosing a coordinate patch
$(U,x^{i}),$ some $\Omega \in \Lambda _{hor}^{k}(T)$ over $\pi ^{-1}(U)$ is
of the form $\Omega _{i_{1}i_{2}...i_{k}}(x,\xi _{x})$ where $\xi _{x}$
denotes an \textit{arbitrary }point in the fiber $\pi ^{-1}(x)$ so that $x$
and $\xi _{x}$ are \textit{independent }variables. For simplicity of
notation, we write $\xi $ for $\xi _{x}.$ We call $\Omega $ smooth if its
components are smooth functions. Henceforth we always assume that our forms
are smooth.

We can express $\Omega $ also as

\begin{equation}
\frac{1}{k!}\sum \text{ }\Omega _{i_{1}i_{2}...i_{k}}(x,\xi )\text{ }%
dx^{i_{1}}\wedge dx^{i_{2}}\wedge ...\wedge dx^{i_{k}}
\end{equation}

Using Einstein summation convention, our shorthand notation for (6) will be $%
\Omega _{I}(x,\xi )dx^{I}$ or just $\Omega _{I}(x,\xi )$ with the obvious
meaning of the index symbol $I.$ We define the total derivative $\widehat{D}%
_{r}$ with respect to the variable $x^{r}$ of a horizontal $0$-form by the
formula

\begin{equation}
\widehat{D}_{r}\Omega (x,\xi )\overset{def}{=}\frac{\partial \Omega (x,\xi )%
}{\partial x^{r}}+\frac{\partial \Omega (x,\xi )}{\partial \xi ^{a}}\widehat{%
\Gamma }_{rb}^{a}(x)\xi (x)^{b}
\end{equation}%
\ \

In other words, we pretend that $\xi $ depends on $x,$ differentiate $\Omega
(x,\xi (x))$ with respect to $x^{r}$ and formally substitute $\frac{\partial
\xi ^{i}}{\partial x^{r}}$ from the equation $\widehat{\nabla }\xi =0.$ Note
that a horizontal form on $J^{\infty }T$ is of the form $\Omega
_{I}(x^{i},\xi ^{i},\xi _{j_{1}}^{i},\xi _{j_{1}j_{2}}^{i},...,\xi
_{j_{1}j_{2}....j_{s}}^{i})$ for arbitrary $s$ but $\widetilde{\nabla }\xi
=0 $ (or $\widehat{\nabla }\xi =0)$ implies that all derivatives of $\xi $
are determined by $\xi .$ Henceforth we sometimes omit the notation for the
variables $x,$ $\xi $ from our formulas and write, for instance, (7) as $%
\widehat{D}_{r}\Omega \overset{def}{=}\frac{\partial \Omega }{\partial x^{r}}%
+\frac{\partial \Omega }{\partial \eta ^{a}}\widehat{\Gamma }_{rb}\xi ^{b}$.

Now we define an operator $\widehat{d}:$ $\Lambda _{hor}^{k}(T)\rightarrow
\Lambda _{hor}^{k+1}(T)$ by the formula

\begin{eqnarray}
&&\left( \widehat{d}\Omega \right) _{ri_{1}i_{2}...i_{k}}\overset{def}{=}%
\left[ \widehat{D}_{r}\Omega _{i_{1}i_{2}...i_{k}}\right]
_{[ri_{1}i_{2}...i_{k}]} \\
&=&\widehat{D}_{r}\Omega _{i_{1}i_{2}...i_{k}}-\widehat{D}_{i_{1}}\Omega
_{ri_{2}...i_{k}}-\widehat{D}_{i_{2}}\Omega _{i_{1}r...i_{k}}...-\widehat{D}%
_{i_{k}}\Omega _{i_{1}i_{2}...r}  \notag
\end{eqnarray}

Equivalently, $\widehat{d}\Omega \overset{def}{=}\widehat{D}_{r}\Omega _{I}$
$dx^{r}\wedge dx^{I}.$

\begin{proposition}
$\widehat{d}\circ \widehat{d}=0$
\end{proposition}

Proof: It suffices to show that $\widehat{D}_{s}\widehat{D}_{r}\Omega $ is
symmetric in $s,r$ for a horizontal $0$-form $\Omega .$ Applying $\widehat{D}%
_{s}$ to (7) gives

\begin{eqnarray}
&&\frac{\partial \Omega }{\partial x^{s}\partial x^{r}}+\frac{\partial
\Omega }{\partial \xi ^{a}\partial x^{r}}\widehat{\Gamma }_{sb}\xi ^{b}+%
\frac{\partial \Omega }{\partial x^{s}\partial \xi ^{a}}\widehat{\Gamma }%
_{rb}\xi ^{b}+\frac{\partial \Omega }{\xi ^{c}\partial \xi ^{a}}\widehat{%
\Gamma }_{rb}\xi ^{b}\widehat{\Gamma }_{sd}^{c}\xi ^{d}  \notag \\
&&+\frac{\partial f}{\partial \xi ^{a}}\left( \frac{\partial \widehat{\Gamma
}_{rb}^{a}}{\partial x^{s}}+\widehat{\Gamma }_{sb}^{c}\widehat{\Gamma }%
_{rc}^{a}\right) \xi ^{b}
\end{eqnarray}

The sum of the first four terms is clearly symmetric in $s,r.$ The last term
is also symmetric in view of $\widehat{\mathfrak{R}}_{sr,b}^{a}=0.$ \ $%
\square $

Thus we get the complex

\begin{equation}
\Lambda _{hor}^{0}(T)\overset{\widehat{d}}{\longrightarrow }\Lambda
_{hor}^{1}(T)\overset{\widehat{d}}{\longrightarrow }\Lambda _{hor}^{2}(T)%
\overset{\widehat{d}}{\longrightarrow }....\overset{\widehat{d}}{%
\longrightarrow }\Lambda _{hor}^{n}(T)
\end{equation}

Note that (10) can be constructed on any parallelizable manifold using $%
\widetilde{D}$ instead of $\widehat{D}$ because we always have $\widetilde{%
\mathfrak{R}}=0.$

Unfortunately, the cohomology groups of (10) turn out to be infinite
dimensional. To see this, we consider now the kernel of the first operator
in (10). First, note that $\widehat{\varepsilon }(p,q)$ (or $\widetilde{%
\varepsilon }(p,q))$ defines an isomorphism of tensor spaces $\widehat{%
\varepsilon }(p,q)_{\ast }:T_{s}^{r}(p)\rightarrow T_{s}^{r}(q).$

\begin{definition}
$\Omega \in \Lambda _{hor}^{0}(T)$ is $\widehat{\varepsilon }$-invariant if $%
\Omega (q,\widehat{\varepsilon }(p,q)_{\ast }\xi )=\Omega (p,\xi ),$ $p,q\in
U,$ $\xi \in \pi ^{-1}(p)$
\end{definition}

Unlike $\widetilde{\varepsilon }$-invariance, $\widehat{\varepsilon }$%
-invariance is a local concept unless $(M,\widetilde{\varepsilon })$ is
globalizable. Equivalently, we may fix $p$ arbitrarily and let $q$ vary in
the condition of Definition 4. Therefore, choosing coordinates around $p,$ $%
q,$ we write this invariance condition as

\begin{equation}
\Omega (x,\eta )=\Omega (p,\xi )\text{ \ \ \ }\eta ^{i}=\widehat{\varepsilon
}_{a}^{i}(p,x)\xi ^{a}
\end{equation}

\QTP{Body Math}
We denote the vector space of $\widehat{\varepsilon }$-invariant $0$-forms
by $\widehat{\Lambda _{hor}^{0}(T)}.$

\begin{proposition}
The kernel of the first operator $\widehat{d}$ in (10) coincides with $%
\widehat{\Lambda _{hor}^{0}(T)}.$
\end{proposition}

Proof: Differentiation of (11) with respect to $x^{r}$ at $x=p$ gives$\left(
\widehat{D}_{r}\Omega \right) (p)=0.$ Since $p$ is arbitrary, we conclude $%
\widehat{D}_{r}\Omega =(\widehat{d}\Omega )_{r}=0.$ Conversely assume $%
\widehat{D}_{r}\Omega =0,$ fix $p,\xi ,x$ and consider the LHS of (11)
defined by the condition in (11). We want to show the equality in (11). Now $%
\widehat{D}_{r}\Omega =0$ implies $\frac{\partial \Omega (x,\eta )}{\partial
x^{r}}=0$ so that $\Omega (x,\eta )$ is independent of $x.$ Setting $x=p$ we
get (11). \ $\square $

We will use the principle in the proof of Proposition 5 several times later
on without giving further details. Now since some $\Omega \in \widehat{%
\Lambda _{hor}^{0}(T)}$ is globally determined by its values on some fiber $%
\pi ^{-1}(p)$ and the vector space of smooth functions on $\pi ^{-1}(p)$ is
infinite dimensional, we conclude $\dim \widehat{\Lambda _{hor}^{0}(T)}%
=\infty .$ This deficiency of (10) forces us to assume the linearity of our
horizontal forms on the fibers. Surprisingly, if $(M,\widetilde{\varepsilon }%
)$ is globalizable and $M$ is compact, this assumption makes the cohomology
of (10) finite dimensional and even computable as we will see shortly.

\begin{definition}
A horizontal $k$-form is linear if it is a linear function on the fibers $%
\pi ^{-1}(x)$ of $\pi :T\rightarrow M.$
\end{definition}

A horizontal linear $k$-form $\omega $ is locally of the form $\omega (x,\xi
)=$ $\omega _{a,i_{1}i_{2}...i_{k}}(x)\xi ^{a}$ where $\omega
_{a,i_{1}i_{2}...i_{k}}$ is a tensor alternating in the indices $%
i_{1}i_{2}...i_{k}$. Therefore a horizontal linear $k$-form on $T\rightarrow
M$ is simply a section of $T^{\ast }(M)\otimes \Lambda ^{k}(M)\rightarrow M,$
whose total space (and also the space of its sections) will be denoted
simply by $T^{\ast }\otimes \Lambda ^{k}.$ In particular, a horizontal
linear $0$-form is an ordinary $1$-form and (7) becomes

\begin{equation}
\left( \widehat{D}_{r}\omega \right) (x,\xi )=\left( \frac{\partial \omega
_{a}}{\partial x^{r}}+\omega _{b}\widehat{\Gamma }_{ra}\right) \xi
^{a}=\left( \widehat{\nabla }_{r}\omega _{a}\right) \xi ^{a}
\end{equation}

Thus we get the subcomplex

\begin{equation}
0\longrightarrow \widehat{\Lambda ^{1}}\longrightarrow \Lambda ^{1}=T^{\ast }%
\overset{\widehat{d}}{\longrightarrow }T^{\ast }\otimes \Lambda ^{1}\overset{%
\widehat{d}}{\longrightarrow }T^{\ast }\otimes \Lambda ^{2}\overset{\widehat{%
d}}{\longrightarrow }....\overset{\widehat{d}}{\longrightarrow }T^{\ast
}\otimes \Lambda ^{n}
\end{equation}%
of (10) and clearly $\dim $ $\widehat{\Lambda ^{1}}=\dim M.$

\begin{definition}
(13) is the horizontal linear complex (LHC) of the local Lie group $(M,%
\widetilde{\varepsilon }).$
\end{definition}

A horizontal $k$-form $\Omega $ (not necessarily linear) is $\widetilde{%
\varepsilon }$-invariant if $\widetilde{\varepsilon }(p,x)_{\ast }\Omega
(x,\eta )=\Omega (p,\widetilde{\varepsilon }(p,x)_{\ast }\xi )$. In
coordinates this condition is

\begin{equation}
\widetilde{\varepsilon }_{i_{1}}^{a_{1}}(p,x)...\widetilde{\varepsilon }%
_{i_{k}}^{a_{k}}(p,x)\Omega _{a_{1}...a_{k}}(x,\eta )=\Omega
_{i_{1}...i_{k}}(p,\xi )\text{ \ \ \ }\eta ^{i}=\widetilde{\varepsilon }%
_{a}^{i}(p,x)\xi ^{a}
\end{equation}

Differentiation of (14) at $x=p$ gives

\begin{equation}
\widetilde{\Gamma }_{r\text{ }i_{1}}^{a_{1}}\Omega
_{a_{1}i_{2}...i_{k}}+....+\widetilde{\Gamma }_{r\text{ }i_{k}}^{a_{k}}%
\Omega _{i_{1}i_{2}...a_{k}}+\frac{\partial \Omega _{i_{1}...i_{k}}}{%
\partial x^{r}}+\frac{\partial \Omega _{i_{1}...i_{k}}}{\partial \eta ^{a}}%
\widetilde{\Gamma }_{rb}^{a}\xi ^{b}=0
\end{equation}

We denote the expression on the LHS of (15) by $\widetilde{\square }%
_{r}\Omega _{i_{1}...i_{k}}$ and call $\widetilde{\square }_{r}$ the $%
\widetilde{}$ -covariant derivative of $\Omega $ with respect to $x^{r}.$
Since $p,\xi $ are arbitrary in (15), if $\Omega $ is $\widetilde{%
\varepsilon }$-invariant then $\widetilde{\square }\Omega =0.$ Converse also
holds and the proof is identical with the proof of Proposition 5. It is
crucial to observe that $\widetilde{D}_{r}\Omega _{I}$ is not a linear
object like $\widetilde{\square }_{r}\Omega _{I}$ unless $\Omega _{I}$ is a $%
0$-form and $\widetilde{D}\Omega =\widetilde{\square }\Omega $ for a $0$%
-form $\Omega .$ This is the reason why we alternate as in (8) to get the
linear object $(\widetilde{d}\Omega )_{ri_{1}i_{2}...i_{k}}$ from $%
\widetilde{D}_{r}\Omega _{i_{1}i_{2}...i_{k}}.$ As another crucial fact, if
we replace, for instance, $\widehat{D}_{r}$ by $\widehat{\square }_{r}$ in
(8), we get a quite different linear object unless $\Omega $ is a $0$-form
but the new operator obtained this way does not give a complex like (10) due
to the presence of torsion.

To summarize, the space of $\widetilde{\varepsilon }$-invariant (horizontal)
linear $k$-forms is $\widetilde{T^{\ast }\otimes \Lambda ^{k}.}$ For $\omega
=\left( \omega _{a,i_{1}i_{2}...i_{k}}\right) \in T^{\ast }\otimes \Lambda
^{k}$ (12) and (15) give

\begin{eqnarray}
\widetilde{\square }_{r}\left( \omega _{b,i_{1}i_{2}...i_{k}}\xi ^{b}\right)
&=&\widetilde{\Gamma }_{r\text{ }i_{1}}^{a_{1}}\omega
_{b,a_{1}i_{2}...i_{k}}\xi ^{b}+....+\widetilde{\Gamma }_{r\text{ }%
i_{k}}^{a_{k}}\omega _{b,i_{1}i_{2}...a_{k}}\xi ^{b} \\
&&+\frac{\partial \omega _{b,i_{1}...i_{k}}}{\partial x^{r}}\xi ^{b}+%
\widetilde{\Gamma }_{rb}^{a}\omega _{a,i_{1}i_{2}...a_{k}}\xi ^{b}  \notag \\
&=&\left( \widetilde{\nabla }_{r}\omega _{a,i_{1}...i_{k}}\right) \xi ^{b}
\notag
\end{eqnarray}

\

So $\widetilde{\square }\omega =\widetilde{\nabla }\omega $ where we
interpret $\omega $ as a horizontal linear $k$-form in $\widetilde{\square }%
\omega $ and as a section of $T^{\ast }\otimes \Lambda ^{k}\rightarrow M$ in
$\widetilde{\nabla }\omega .$ Now if $\widetilde{\square }\omega =0,$ then
(16) gives

\begin{equation}
\widehat{D}_{r}\omega _{a,i_{1}i_{2}...a_{k}}=-\widetilde{\Gamma }_{r\text{ }%
i_{1}}^{a_{1}}\omega _{b,a_{1}i_{2}...i_{k}}-....-\widetilde{\Gamma }_{r%
\text{ }i_{k}}^{a_{k}}\omega _{b,i_{1}i_{2}...a_{k}}-\omega
_{a,i_{1}i_{2}...a_{k}}\widetilde{T}_{rb}^{a}
\end{equation}

The proof of the next Proposition is almost identical with the proof of
Proposition 7 in [13]. This is not surprising for if we replace $T^{\ast }$
in (13) by $T$, then (13) becomes (30) in [13].

\begin{proposition}
If $\widetilde{\square }\omega =0$, then $\widetilde{\square }\widehat{d}%
\omega =0.$ Therefore $\widehat{d}:\widetilde{T^{\ast }\otimes \Lambda ^{k}}%
\rightarrow $ $\widetilde{T^{\ast }\otimes \Lambda ^{k+1}}$
\end{proposition}

Proof (sketch): We observe that each term in the alternation $\left[ \omega
_{a,i_{1}i_{2}...a_{k}}\widetilde{T}_{rb}^{a}\right] _{[ri_{1}...i_{k}]}$ of
the last term of (17) is a tensor. Applying $\widetilde{\nabla }_{s}$ to
each such term gives zero by (2). The alternation $\left[ \widetilde{\Gamma }%
_{r\text{ }i_{1}}^{a_{1}}\omega _{b,a_{1}i_{2}...i_{k}}+....+\widetilde{%
\Gamma }_{r\text{ }i_{k}}^{a_{k}}\omega _{b,i_{1}i_{2}...a_{k}}\right]
_{[ri_{1}...i_{k}]}$ is a sum of terms of the form $\widetilde{T}_{\ast \ast
}^{a}\omega _{b,\ast ...a....\ast }$ and we argue as before. \ $\square $

Let $\widetilde{\widehat{\Lambda ^{1}}}$ denote $\widehat{\varepsilon }$%
-invariant $1$-forms which are also $\widetilde{\varepsilon }$-invariant.
Clearly $\widetilde{\widehat{\Lambda ^{1}}}=\widehat{\widetilde{\Lambda ^{1}}%
}=\widehat{\Lambda ^{1}}\cap \widetilde{\Lambda ^{1}}.$ Now Proposition 8
gives the subcomplex

\begin{equation}
0\longrightarrow \widetilde{\widehat{\Lambda ^{1}}}\longrightarrow
\widetilde{\Lambda ^{1}}\overset{\widehat{d}}{\longrightarrow }\widetilde{%
T^{\ast }\otimes \Lambda ^{1}}\overset{\widehat{d}}{\longrightarrow }....%
\overset{\widehat{d}}{\longrightarrow }\widetilde{T^{\ast }\otimes \Lambda
^{n}}
\end{equation}%
of (13) which localizes at any point $p\in M$ and can therefore reduces to
algebra.

\begin{definition}
(18) is the invariant linear horizontal complex (ILHC) of $(M,\widetilde{%
\varepsilon }).$
\end{definition}

If $(M,\widetilde{\varepsilon })$ is globalizable and $M$ compact, then $M$
admits a measure invariant under both the global transformation groups $%
\widetilde{\mathcal{G}}$ and $\widehat{\mathcal{G}}$ and the standard
averaging process over $M$ proves that the inclusion of (18) in (13) induces
isomorphism in cohomology in positive degrees. However, (18) computes the
cohomology of the Lie algebra $\mathfrak{g}\overset{def}{=}\widetilde{\Theta
}$ with coefficients $\mathfrak{g}^{\ast }=\widetilde{\Lambda ^{1}}$ in the
same way as (30) in [13] computes the cohomology of $\mathfrak{g}$ with
coefficients $\mathfrak{g}.$ Thus we conclude

\begin{proposition}
If the local Lie group $(M,\widetilde{\varepsilon })$ is globalizable and $M$
is compact, then the $k$' th cohomology groups of (18) and (13) are both
isomorphic to $H^{k}(\mathfrak{g,g}^{\ast })$ in positive degrees where $%
\mathfrak{g}\overset{def}{=}\widetilde{\Theta }$ and $\mathfrak{g}^{\ast }=%
\widetilde{\Lambda ^{1}}.$
\end{proposition}

Finally, we remark that our construction in this section works if we replace
$T^{\ast }\rightarrow M$ by the $(r,s)$-tensor bundle $T_{s}^{r}\rightarrow
M $ and in fact by any natural vector bundle $E\rightarrow M$ of order one
(see [10] and the references there for natural bundles). For $%
T_{s}^{r}\rightarrow M,$ (18) computes $H^{\ast }(\mathfrak{g,}T_{s}^{r}(%
\mathfrak{g})).$

\section{The nonlinear horizontal complex}

\begin{definition}
A nonlinear horizontal $k$-form $\omega $ on $M\times M$ assigns to $%
(p,q)\in M\times M$ an ordinary $k$-form at $p.$
\end{definition}

Note that $\omega $ can be defined also as a function $\omega :M\times
M\rightarrow \widetilde{\Lambda ^{k}(M)\text{ }}$ since elements of $%
\widetilde{\Lambda ^{k}(M)\text{ }}$ are globally determined by their values
at any point. We denote the space of (nonlinear) horizontal $k$-forms by $%
\Lambda _{hor}^{k}(M\times M).$ Choosing coordinates $p\in (U,x^{i}),$ $q\in
(V,y^{i}),$ we write $\omega $ as $\omega _{i_{1}i_{2}...i_{k}}(x,y).$ There
is an ambiguity with this notation: it does not specify the coordinates to
which the $k$-form indices $i_{1},i_{2},$ $...,i_{k}$ refer to. Except in
the proof of Proposition 22, we agree that they refer to the coordinates
around the source point $p.$ Note that a choice of coordinates around some
point canonically defines coordinates around all points if $(M,\widetilde{%
\varepsilon })$ is a local Lie group.

In view of (3), we define the total differentiation operator $\widetilde{D}%
:\Lambda _{hor}^{0}(M\times M)\rightarrow \Lambda _{hor}^{1}(M\times M)$ by

\begin{equation}
(\widetilde{D}\theta )_{r}(x,y)\overset{def}{=}\frac{\partial \theta (x,y)}{%
\partial x^{r}}+\frac{\partial \theta (x,y)}{\partial y^{a}}\widetilde{%
\varepsilon }_{r}^{a}(x,y)
\end{equation}

Now we define $\widetilde{d}:\Lambda _{hor}^{k}(M\times M)\rightarrow
\Lambda _{hor}^{k+1}(M\times M)$ by

\begin{equation}
\widetilde{d}\omega \overset{def}{=}\left( \widetilde{D}\omega _{I}\right)
_{r}\text{ }dx^{r}\wedge dx^{I}
\end{equation}

\begin{proposition}
If $(M,\varepsilon )$ is a local Lie group, then $\widetilde{d}\circ
\widetilde{d}=0$
\end{proposition}

Proof: Writing $(\widetilde{D}f)_{r}$ as $\widetilde{D}_{r}f$ and applying $%
\widetilde{D}_{s}$ to (19) gives

\begin{eqnarray}
&&\frac{\partial ^{2}\theta (x,y)}{\partial x^{s}\partial x^{r}}+\frac{%
\partial ^{2}\theta (x,y)}{\partial y^{a}\partial x^{r}}\varepsilon
_{s}^{a}(x,y)+\frac{\partial \theta (x,y)}{\partial x^{s}\partial y^{a}}%
\varepsilon _{r}^{a}(x,y) \\
&&+\frac{\partial \theta (x,y)}{\partial y^{b}\partial y^{a}}\varepsilon
_{s}^{a}(x,y)\varepsilon _{r}^{a}(x,y)+\frac{\partial \theta (x,y)}{\partial
y^{a}}\left( \frac{\partial \varepsilon _{r}^{a}(x,y)}{\partial x^{s}}+\frac{%
\partial \varepsilon _{r}^{a}(x,y)}{\partial y^{b}}\varepsilon
_{s}^{b}(x,y)\right)  \notag
\end{eqnarray}%
which is symmetric in $s,r$ since $\widetilde{\mathcal{R}}_{sr}^{a}(x,y)=0.$
\ $\square $

Thus we get the complex
\begin{equation}
\Lambda _{hor}^{0}(M\times M)\overset{\widetilde{d}}{\longrightarrow }%
\Lambda _{hor}^{1}(M\times M)\overset{\widetilde{d}}{\longrightarrow }....%
\overset{\widetilde{d}}{\longrightarrow }\Lambda _{hor}^{n}(M\times M)
\end{equation}

\begin{definition}
(22) is the nonlinear horizontal complex (NHC) of the local Lie group $(M,%
\widetilde{\varepsilon }).$
\end{definition}

The construction of (22) needs the parallelizable manifold $(M,\widetilde{%
\varepsilon })$ together with the assumption $\widetilde{\mathcal{R}}=0.$
The restriction $(U,\widetilde{\varepsilon }_{\mid U})$ of $(M,\widetilde{%
\varepsilon })$ to some $U\subset M$ satisfies both these conditions. Hence
we can meaningfully speak of the restriction of (22) to $U.$ Now the
following question arises naturally

$\mathbf{Q:}$ Is (22) locally exact?

Let $f\in \widetilde{\mathcal{G}}$ be the unique local solution of (3) with
the initial condition $f(p)=q,$ $p\in (U,x^{i})$ and $\left\{ (x,f(x)),\text{
}x\in U\right\} $ be the local graph of $f.$ For $\theta \in \Lambda
_{hor}^{0}(M\times M)$ we consider the restriction $\theta (x,f(x))$ of $%
\theta $ to the graph of $f\in \widetilde{\mathcal{G}}.$

The proof of the next proposition follows easily from the definitions.

\begin{proposition}
The following are equivalent

i) $\theta \in $ $\Lambda _{hor}^{0}(M\times M)$ belongs to the kernel of
the first operator in (22)

ii) The restriction of $\theta $ to the graph of $f$ is constant for all $%
f\in \widetilde{\mathcal{G}}$.
\end{proposition}

To understand this kernel better, it is useful to assume that $(M,\widetilde{%
\varepsilon })$ is globalizable so that $\widetilde{\mathcal{G}}$ is a
global transformation group of $M$ which acts simply transitively. So we may
identify $f\in \widetilde{\mathcal{G}}$ with its graph $\left\{ (p,f(p)),%
\text{ }p\in M\right\} \subset M\times M$. Since $\theta $ is constant on
this graph (we always assume that $M$ is connected), we interpret this
constant value as the value of $\theta $ on $f.$ This identifies the kernel
with the functions $\theta :\widetilde{\mathcal{G}}\rightarrow \mathbb{R}.$

We recall that $g(a,b,x)$ is the unique solution of (3) in the variable $x$
satisfying the initial condition $a\rightarrow b.$ Let $\theta \in \Lambda
_{hor}^{0}(M\times M)$. We call $\theta $ $\widehat{\varepsilon }$-invariant
if $\theta (x,g(p,x,q))=$ $\theta (p,q)$ for $p,q,x\in M.$ Since $g(a,b,x)$
is defined for sufficiently close $a,x$ and $p,q$ are arbitrary in our
definition, we make the flat assumption of globalizability henceforth so
that $\widehat{\mathcal{G}}$ is another global transformation group of $M$
which acts simply transitively. Differentiating $\theta (x,g(p,x,q))=$ $%
\theta (p,q)$ with respect to $x$ at $x=p$ gives

\begin{equation}
\widehat{D}_{r}(p,q)=\frac{\partial \theta }{\partial x^{r}}(p,q)+\frac{%
\partial \theta }{\partial y^{a}}(p,q)\widehat{\varepsilon }(p,q)=0
\end{equation}

Since $p,q$ are arbitrary in (23), we deduce $\widehat{D}\theta =0$ and we
easily show as before that conversely $\widehat{D}\theta =0$ implies the $%
\widehat{\varepsilon }$- invariance of $\theta .$ Let $\widehat{\Lambda
_{hor}^{0}(M\times M)}$ denote the space of $\widehat{\varepsilon }$%
-invariant functions and let $\theta \in \widehat{\Lambda _{hor}^{0}(M\times
M)}$ satisfy $\widetilde{d}\theta =0.$ The proof of Proposition 14 shows
that $\theta $ is constant on the graphs of $h\in \widehat{\mathcal{G}}$ and
therefore may be interpreted as a function $\theta :\widehat{\mathcal{G}}%
\rightarrow \mathbb{R}$. Now fix some $f\in \widetilde{\mathcal{G}}$, some $%
p\in M$ and suppose $f(p)=q.$ Let $x\in M.$ There is a unique $k_{x}\in
\widetilde{\mathcal{G}}$ with $k_{x}(p)=x.$ Therefore $g(p,x,q)=k_{x}(q)=$ $%
\left( k_{x}\circ f\circ k_{x}^{-1}\right) (x)$ and $h$ in (5) is the
transformation $x\rightarrow k_{x}\circ f\circ k_{x}^{-1}(x).$ Hence we
conclude that this transformation belongs to $\widehat{\mathcal{G}}$. Now
since $\theta \in \widehat{\Lambda _{hor}^{0}(M\times M)},$ $\theta
(k_{x}\circ f\circ k_{x}^{-1})$ has the same value $\theta (p,q)$
independent of $x.$ However, since also $\widetilde{d}\theta =0$, we have $%
\theta (p,q)=\theta (f).$ Whence

\begin{equation}
\theta (f)=\theta (k\circ f\circ k^{-1})\text{ \ }f,k\in \widetilde{\mathcal{%
G}}
\end{equation}

Recall that a function on a Lie group which is constant on the conjugacy
classes is called a character function. The trace of a representation is a
character function and these functions play a fundamental role in
representation theory. Let $\mathcal{C(}\widetilde{\mathcal{G}})$ denote the
space of character functions defined by (24). Thus we proved

\begin{proposition}
The sequence
\end{proposition}

\begin{equation}
0\rightarrow \mathcal{C(}\widetilde{\mathcal{G}})=\widetilde{\widehat{%
\Lambda _{hor}^{0}(M\times M)}}\longrightarrow \widehat{\Lambda
_{hor}^{0}(M\times M)}\overset{\widetilde{d}}{\longrightarrow }\Lambda
_{hor}^{1}(M\times M)
\end{equation}

is exact.

\begin{definition}
$\omega \in \Lambda _{hor}^{k}(M\times M)$ is $\widehat{\varepsilon }$%
-invariant if $\widehat{\varepsilon }(p,x)_{\ast }\omega (x,g(p,x,q))=$ $%
\omega (p,q).$
\end{definition}

In coordinates the condition of invariance is

\begin{equation}
\widehat{\varepsilon }_{i_{1}}^{a_{1}}(p,x)...\widehat{\varepsilon }%
_{i_{k}}^{a_{k}}(p,x)\omega _{a_{1}...a_{k}}(x,g(p,x,q))=\omega
_{i_{1}...i_{k}}(p,q)
\end{equation}

Differentiation of (26) with respect to $x^{r}$ at $x=p$ gives (omitting $%
p,q $ from our notation and using the same symbol $\widehat{\square }$ as
before) $\widehat{\square }_{r}\omega _{a_{1}i_{2}...i_{k}}=0$ where

\begin{eqnarray}
&&\widehat{\square }_{r}\omega _{a_{1}i_{2}...i_{k}}\overset{def}{=}\widehat{%
\Gamma }_{ri_{1}}^{a_{1}}\omega _{a_{1}i_{2}...i_{k}}+\widehat{\Gamma }%
_{ri_{2}}^{a_{2}}\omega _{i_{1}a_{2}...i_{k}}+...+\widehat{\Gamma }%
_{ri_{k}}^{a_{k}}\omega _{i_{1}i_{2}...a_{k}}  \notag \\
&&+\frac{\partial \omega _{i_{1}i_{2}...a_{k}}}{\partial x^{r}}+\frac{%
\partial \omega _{i_{1}i_{2}...a_{k}}}{\partial y^{a}}\widehat{\varepsilon }%
_{r}^{a} \\
&=&\widehat{\Gamma }_{ri_{1}}^{a_{1}}\omega _{a_{1}i_{2}...i_{k}}+\widehat{%
\Gamma }_{ri_{2}}^{a_{2}}\omega _{i_{1}a_{2}...i_{k}}+...+\widehat{\Gamma }%
_{ri_{k}}^{a_{k}}\omega _{i_{1}i_{2}...a_{k}}+\widehat{D}_{r}\omega
_{i_{1}i_{2}...i_{k}}  \notag
\end{eqnarray}

As before, some $\omega \in \Lambda _{hor}^{k}(M\times M)$ is $\widehat{%
\varepsilon }$-invariant if and only if $\widehat{\square }\omega =0.$ Using
(27), an argument similar to the proof of Proposition 8 gives

\begin{proposition}
We have the complex
\end{proposition}

\begin{equation}
0\rightarrow \mathcal{C(}\widetilde{\mathcal{G}})\rightarrow \widehat{%
\Lambda _{hor}^{0}(M\times M)}\overset{\widetilde{d}}{\longrightarrow }%
\widehat{\Lambda _{hor}^{1}(M\times M)}\overset{\widetilde{d}}{%
\longrightarrow }....\overset{\widetilde{d}}{\longrightarrow }\widehat{%
\Lambda _{hor}^{n}(M\times M)}
\end{equation}

\begin{definition}
(28) is the invariant nonlinear horizontal complex (INHC) of the local Lie
group $(M,\widetilde{\varepsilon }).$
\end{definition}

Observe that (28) does not restrict to $U\subset M$ since $(U,\widetilde{%
\varepsilon }_{\mid U})$ need not be globalizable even if $(M,\widetilde{%
\varepsilon })$ is. We are unable to express the cohomology of (28) in
positive degrees in terms of some known cohomology groups. We also do not
know any sufficient condition which makes the cohomologies of (28) and (22)
isomorphic. However it is worthwhile to note that elements of $\Lambda
_{hor}^{n}(M\times M)$ may be viewed as functionals on the diffeomorphism
group $Diff(M)$ of $M$ if $M$ is compact. Indeed, if $\omega \in $ $\Lambda
_{hor}^{n}(M\times M)$ and $f\in Diff(M)$, then $\omega (x,f(x))$ defines a
volume form as $x$ ranges over $M$ and therefore can be integrated over $M$
giving the functional $\omega :f\rightarrow \int_{M}\omega (x,f(x)).$ This
suggests to continue (22) one step to the right by the Euler-Lagrange
operator $EL$ but we will not enter this issue here.

\section{The linearization map}

Our purpose in this section is to define a chain map from (22) to (13) which
restricts to the invariant subcomplexes (28), (18).

Let $\omega =\omega _{I}(x,y)$ be a nonlinear horizontal $k$-form and $\xi $
be a tangent vector at $x.$ The idea is to let $y$ approach $x$ along the
tangent vector $\xi $. Since $\omega _{I}(x,x+t\xi )$ and $\omega _{I}(x,x)$
are two ordinary $k$-forms at the \textit{same }point $x,$ $\left[ \frac{%
d\omega _{I}(x,x+t\xi )}{dt}\right] _{t=0}$ is well defined and is an
ordinary $k$-form at $x$ which depends on $\xi $, that is, an element of $%
\Lambda _{hor}^{k}(T).$ It depends linearly on $\xi $ because

\begin{eqnarray}
\left[ \frac{d\omega _{I}(x,x+t\xi )}{dt}\right] _{t=0} &=&\left[ \frac{%
\partial \omega _{I}(x,y)}{\partial y^{a}}\right] _{y=x}\xi ^{a} \\
&=&\frac{\partial \omega _{I}(x,x)}{\partial y^{a}}\xi ^{a}  \notag
\end{eqnarray}

Now (29) defines a map $L:$ $\Lambda _{hor}^{k}(M\times M)\rightarrow
T^{\ast }\otimes \Lambda ^{k}.$

\begin{proposition}
The following diagram commutes
\end{proposition}

\begin{equation}
\begin{array}{ccccccc}
\Lambda _{hor}^{0}(M\times M) & \overset{\widetilde{d}}{\longrightarrow } &
\Lambda _{hor}^{1}(M\times M) & \overset{\widetilde{d}}{\longrightarrow } &
... & \overset{\widetilde{d}}{\longrightarrow } & \Lambda _{hor}^{n}(M\times
M) \\
\downarrow _{L} &  & \downarrow _{L} &  & \downarrow _{L} &  & \downarrow
_{L} \\
T^{\ast }=\Lambda ^{1} & \overset{\widehat{d}}{\longrightarrow } & T^{\ast
}\otimes \Lambda ^{1} & \overset{\widehat{d}}{\longrightarrow } & ... &
\overset{\widehat{d}}{\longrightarrow } & T^{\ast }\otimes \Lambda ^{n}%
\end{array}%
\end{equation}

Proof: Let $\omega \in \Lambda _{hor}^{k}(M\times M).$ To compute $\left(
\widehat{D}_{r}\circ L\right) \omega ,$ we apply $\widehat{D}_{r}$ to (29)
which gives

\begin{equation}
\frac{\partial ^{2}\omega _{I}(x,x)}{\partial x^{r}\partial y^{a}}\xi ^{a}+%
\frac{\partial ^{2}\omega _{I}(x,x)}{\partial y^{r}\partial y^{a}}\xi ^{a}+%
\frac{\partial \omega _{I}(x,x)}{\partial y^{a}}\widehat{\Gamma }%
_{rb}^{a}\xi ^{b}
\end{equation}

To compute $\left( L\circ \widetilde{D}_{r}\right) \omega $ we apply $L$ to
(19) which gives

\begin{equation}
\frac{d}{dt}\left[ \frac{\partial \omega _{I}(x,x+t\xi )}{\partial x^{r}}+%
\frac{\partial \omega _{I}(x,x+t\xi )}{\partial y^{a}}\widetilde{\varepsilon
}_{r}^{a}(x,x+t\xi )\right] _{t=0}
\end{equation}%
and (1) shows that (31) and (32) are equal. \ $\square $

Unfortunately, we do \textit{not }have $L:\widehat{\Lambda
_{hor}^{k}(M\times M)}\longrightarrow \widetilde{T^{\ast }\otimes \Lambda
^{k}}.$ For instance, let $k=2$ and $\omega \in \widehat{\Lambda
_{hor}^{2}(M\times M)}.$ By (27), the condition $\widehat{\square }%
_{r}\omega _{ij}=0$ is

\begin{equation}
\widehat{\Gamma }_{ri}^{a}(x)\omega _{aj}(x,y)+\widehat{\Gamma }%
_{rj}^{a}(x)\omega _{ia}(x,y)+\frac{\partial \omega _{ij}(x,y)}{\partial
x^{r}}+\frac{\partial \omega _{ij}(x,y)}{\partial y^{a}}\widehat{\varepsilon
}_{r}^{a}(x,y)=0
\end{equation}

Setting $y=x+t\xi $ in (33) and differentiating at $t=0$ gives

\begin{eqnarray}
0 &=&\widehat{\Gamma }_{ri}^{a}(x)\frac{\partial \omega _{aj}(x,x)}{\partial
y^{b}}\xi ^{b}+\widehat{\Gamma }_{rj}^{a}(x)\frac{\partial \omega _{ia}(x,x)%
}{\partial y^{b}}\xi ^{b}+\frac{\partial \omega _{ij}(x,x)}{\partial
y^{a}\partial x^{r}}\xi ^{a}  \notag \\
&&+\frac{\partial \omega _{ij}(x,x)}{\partial y^{a}\partial y^{r}}\xi ^{a}+%
\frac{\partial \omega _{ij}(x,x)}{\partial y^{a}}\widehat{\Gamma }%
_{rb}^{a}(x)\xi ^{b}  \notag \\
&=&\widehat{\Gamma }_{ri}^{a}(x)\frac{\partial \omega _{aj}(x,x)}{\partial
y^{b}}\xi ^{b}+\widehat{\Gamma }_{rj}^{a}(x)\frac{\partial \omega _{ia}(x,x)%
}{\partial y^{b}}\xi ^{b}+\frac{\partial }{\partial x^{r}}\left[ \frac{%
\partial \omega _{ij}(x,x)}{\partial y^{a}}\right] \xi ^{a}  \notag \\
&&+\frac{\partial \omega _{ij}(x,x)}{\partial y^{a}}\widehat{\Gamma }%
_{rb}^{a}(x)\xi ^{b}  \notag \\
&=&\widehat{\Gamma }_{ri}^{a}(x)\left( L\omega \right) _{aj}+\widehat{\Gamma
}_{rj}^{a}(x)\left( L\omega \right) _{ia}+\widehat{D}_{r}\left( L\omega
\right) _{ij}  \notag \\
&=&\widehat{\square }_{r}\left( L\omega \right) _{ij}
\end{eqnarray}%
whereas what we want is $\widetilde{\square }_{r}\left( L\omega \right)
_{ij}=0.$ Clearly we can replace $\widehat{}$ with $\widetilde{}$ in (34).
This makes it necessary to consider forms which are both $\widetilde{%
\varepsilon }$ and $\widehat{\varepsilon }$ invariant. Now the proof of
Proposition 8 shows that other than (18) we also have the subcomplex

\begin{equation}
\widehat{\Lambda ^{1}}\overset{\widehat{d}}{\longrightarrow }\widehat{%
T^{\ast }\otimes \Lambda ^{1}}\overset{\widehat{d}}{\longrightarrow }....%
\overset{\widehat{d}}{\longrightarrow }\widehat{T^{\ast }\otimes \Lambda ^{n}%
}
\end{equation}%
of (13). Observe that the first operator in (35) vanishes on $\widehat{%
\Lambda ^{1}}$. The interpretations of (18) and (35) in the modern formalism
are somewhat intriguing: (35) computes the cohomology of $\mathfrak{g=}%
\widehat{\Theta }$ with coefficients $\widehat{\Lambda ^{1}}$ but the
\textit{representation is trivial. }So (35) computes $n$-copies of the
cohomology of $\widehat{\Theta }$ with trivial coefficients $\mathbb{R}$.
However the representation in (18) is \textquotedblleft
honest\textquotedblright\ (which comes, of course, from the Lie derivative $%
\mathcal{L}_{\widetilde{\xi }}$, see (49) below). Now (18) and (35) give the
subcomplex

\begin{equation}
\widetilde{\widehat{\Lambda ^{1}}}\overset{\widehat{d}}{\longrightarrow }%
\widetilde{\widehat{T^{\ast }\otimes \Lambda ^{1}}}\overset{\widehat{d}}{%
\longrightarrow }....\overset{\widehat{d}}{\longrightarrow }\widetilde{%
\widehat{T^{\ast }\otimes \Lambda ^{n}}}
\end{equation}%
(note that a biinvariant form need not be closed in the presence of a
\textquotedblleft representation\textquotedblright\ as can easily be seen
from the second formula in (16) which the reader may compare to 19 in [4]).
We call (36) the biinvariant subcomplex. Similarly we construct the
biinvariant nonlinear complex. The computation in (34) now gives the
following

\begin{proposition}
$L$ restricts to biinvariant subcomplexes.
\end{proposition}

\section{Double complexes}

The main idea of the nonlinear double complex is quite simple: we let the
number of copies of $M$ be arbitrary, modify $\widetilde{d}$ accordingly and
define the vertical operator $\delta $ by the well known formula from
topology and group cohomology.\ Some formulas look quite complicated in
coordinates even though they are straightforward generalizations of our
previous formulas and state some facts which are evident at this stage. For
this reason our treatment will be short.

\begin{definition}
A nonlinear horizontal $k$-form $\omega $ on $M\times M\times ...\times M$ $%
(m$ copies, $m\geq 1)$ assigns to $(p,q,...,t)\in M\times M\times ...\times
M $ an ordinary $k$-form at $p.$
\end{definition}

Let $M^{(m)}$ denote $M\times M\times ...\times M$ $(m$ copies, $m\geq 1).$
So $\omega $ is a function $\omega :M^{(m)}\rightarrow \widetilde{\Lambda
^{k}(M)}.$ According to the formalism of groupoids, the groupoid $\widetilde{%
\varepsilon }(M\times M)\subset J_{1}(M\times M)$ has a representation on
the vector bundle $\Lambda ^{k}(M)\rightarrow M$ and a (nonlinear)
horizontal $k$-form is an $(m-1)$-composable string. Thus we can define the
differentiable cohomology of this groupoid with coefficients $\Lambda
^{k}(M) $ (see [5] for details).

We denote the space of horizontal $k$-forms on $M^{(m)}$ by $\Lambda
_{hor}^{m,k}(M^{(m)}).$Choosing coordinates $%
(U,x^{i}),(V,y^{i}),...,(W,z^{i})$ around $p,q,...,t,$ we express $\omega $
as $\omega _{I}(x,y,...,z)$ where the index $I$ refers to $(x^{i})$ as
before. We define

\begin{eqnarray}
&&\widetilde{D}_{r}\omega _{i_{1}...i_{k}}(x,y,...,z)\overset{def}{=}\frac{%
\partial \omega _{i_{1}...i_{k}}(x,y,...,z)}{\partial x^{r}}+ \\
&&\frac{\partial \omega _{i_{1}...i_{k}}(x,y,...,z)}{\partial y^{a}}%
\widetilde{\varepsilon }_{r}^{a}(x,y)+....+\frac{\partial \omega
_{i_{1}...i_{k}}(x,y,...,z)}{\partial z^{a}}\widetilde{\varepsilon }%
_{r}^{a}(x,z)  \notag
\end{eqnarray}

Using $\widetilde{D}$ we now define $\widetilde{d}:\Lambda
_{hor}^{m,k}(M^{(m)})\rightarrow \Lambda _{hor}^{m,k+1}(M^{(m)}).$ Since $%
\widetilde{\mathcal{R}}=0$ identically on $M\times M$, we get the complex

\begin{equation}
\Lambda _{hor}^{m,0}(M^{(m)})\overset{\widetilde{d}}{\longrightarrow }%
\Lambda _{hor}^{m,1}(M^{(m)})\overset{\widetilde{d}}{\longrightarrow }....%
\overset{\widetilde{d}}{\longrightarrow }\Lambda _{hor}^{m,n}(M^{(m)})
\end{equation}

For $m=1,$ (38) reduces to the ordinary de Rham complex on $M$ and for $m=2$
it reduces to (22). If $(M,\widetilde{\varepsilon })$ is globalizable, which
we assume henceforth, the kernel of the first operator in (38) can be
identified with functions $\theta :\widetilde{\mathcal{G}}%
^{(m-1)}\rightarrow \mathbb{R}$ where we set $\widetilde{\mathcal{G}}^{(0)}=%
\mathbb{R}$.

Now we define $\delta :\Lambda _{hor}^{m,k}\left( M^{(m)}\right) \rightarrow
\Lambda _{hor}^{m+1,k}\left( M^{(m)}\right) $ by the well known formula

\begin{equation}
(\delta \omega )(p_{0},p_{1},...,p_{m})\overset{def}{=}\sum_{0\leq i\leq
m}(-1)^{i}\omega (p_{0},p_{1},...,\overset{(o)}{p_{i}},...,p_{m})
\end{equation}%
where $\overset{(o)}{p_{i}}$ indicates the omission of the point $p_{i}.$ In
(39) we interpret $\omega $ as a function $\omega :M^{(m)}\rightarrow
\widetilde{\Lambda ^{k}(M)}.$ Clearly $\delta ^{2}=0.$ Thus we get the
diagram%
\begin{equation}
\begin{tabular}{lllllll}
&  &  &  &  &  &  \\
$\uparrow _{\delta }$ &  & $\uparrow _{\delta }$ &  & $\uparrow _{\delta }$
&  & $\uparrow _{\delta }$ \\
$....$ & $\overset{\widetilde{d}}{\longrightarrow }$ & $....$ & $\overset{%
\widetilde{d}}{\longrightarrow }$ & $....$ & $\overset{\widetilde{d}}{%
\longrightarrow }$ & $.....$ \\
$\uparrow _{\delta }$ &  & $\uparrow _{\delta }$ &  & $\uparrow _{\delta }$
&  & $\uparrow _{\delta }$ \\
$\Lambda _{hor}^{3,0}(M^{(3)})$ & $\overset{\widetilde{d}}{\longrightarrow }$
& $\Lambda _{hor}^{3,1}(M^{(3)})$ & $\overset{\widetilde{d}}{\longrightarrow
}$ & $.....$ & $\overset{\widetilde{d}}{\longrightarrow }$ & $\Lambda
_{hor}^{3,n}(M^{(3)})$ \\
$\uparrow _{\delta }$ &  & $\uparrow _{\delta }$ &  & $\uparrow _{\delta }$
&  & $\uparrow _{\delta }$ \\
$\Lambda _{hor}^{2,0}(M^{(2)})$ & $\overset{\widetilde{d}}{\longrightarrow }$
& $\Lambda _{hor}^{2,1}(M^{(2)})$ & $\overset{\widetilde{d}}{\longrightarrow
}$ & $.....$ & $\overset{\widetilde{d}}{\longrightarrow }$ & $\Lambda
_{hor}^{2,n}(M^{(2)})$ \\
$\uparrow _{\delta }$ &  & $_{{}}\uparrow _{\delta }$ &  & $\uparrow
_{\delta }$ &  & $\uparrow _{\delta }$ \\
$\Lambda _{hor}^{1,0}(M)$ & $\overset{\widetilde{d}}{\longrightarrow }$ & $%
\Lambda _{hor}^{1,1}(M)$ & $\overset{\widetilde{d}}{\longrightarrow }$ & $%
.....$ & $\overset{\widetilde{d}}{\longrightarrow }$ & $\Lambda
_{hor}^{1,n}(M)$%
\end{tabular}%
\end{equation}

\begin{proposition}
The diagram (40) commutes
\end{proposition}

Proof: We check the commutativity of the square

\begin{equation}
\begin{tabular}{lll}
$\Lambda _{hor}^{3,k}(M^{3})$ & $\overset{\widetilde{d}}{\longrightarrow }$
& $\Lambda _{hor}^{3,k}(M^{3})$ \\
$\uparrow \delta $ &  & $\uparrow \delta $ \\
$\Lambda _{hor}^{2,k}(M)$ & $\overset{\widetilde{d}}{\longrightarrow }$ & $%
\Lambda _{hor}^{2,k+1}(M)$%
\end{tabular}%
\end{equation}%
and the general case is similar. For $\omega _{I}\in \Lambda
_{hor}^{2,k}(M), $ we have

\begin{equation}
(\delta \omega )_{I}(x,y,z)=\omega _{I}(y,z)-\omega _{I}(x,z)+\omega
_{I}(x,y)
\end{equation}

We should be careful with (42): $I$ refers to $(x^{i})$ and $\omega
_{I}(y,z) $ denotes the value of $\omega (y,z)\in \widetilde{\Lambda ^{k}(M)}
$ at $x.$ Now we assume $y=y(x)$ and $z=z(y)=z(y(x))$ belong to $\widetilde{%
\mathcal{G}}$ with $y(p)=q,$ $z(q)=o$, substitute $y(x)$, $z(x)$ into (42)
and differentiate (42) with respect to $x^{r}$ at $x=p.$ The result is

\begin{eqnarray}
\widetilde{D}_{r}(\delta \omega )_{I}(p,q,o) &=&(\widetilde{D}_{r}\omega
_{I})(q,o)-(\widetilde{D}_{r}\omega _{I})(p,o)+(\widetilde{D}_{r}\omega
_{I})(p,q)  \notag \\
&=&\left( \delta \left( \widetilde{D}_{r}\omega _{I}\right) \right) (p,q,o)
\end{eqnarray}%
and (43) implies $\widetilde{d}\circ \delta =\delta \circ $ $\widetilde{d}.$
\ $\square $

\begin{definition}
The diagram (40) is the nonlinear horizontal double complex (NHDC) of the
local Lie group.
\end{definition}

Some $\omega \in \Lambda _{hor}^{m,k}(M^{(m)})$ is $\widehat{\varepsilon }$%
-invariant if $\widehat{\varepsilon }(p,x)_{\ast }\omega
(x,g(p,x,q),...,g(p,x,t))=\omega (p,q,...,t).$ This condition is (26) in
coordinates except that we should take also the other components into
account. Differentiation of this formula at $x=p$ gives (27) where $\widehat{%
D}_{r}\omega _{i_{1}i_{2}...i_{k}}$ is defined by (35). In this way we get
the subcomplex

\begin{equation}
\widehat{\Lambda _{hor}^{m,0}(M^{m})}\overset{\widetilde{d}}{\longrightarrow
}\widehat{\Lambda _{hor}^{m,1}(M^{m})}\overset{\widetilde{d}}{%
\longrightarrow }....\overset{\widetilde{d}}{\longrightarrow }\widehat{%
\Lambda _{hor}^{m,n}(M^{m})}
\end{equation}

For $m=1,$ $\widehat{\Lambda _{hor}^{k,1}(M^{1})}=\widehat{\Lambda ^{k}}$, $%
\widetilde{d}$ is the ordinary exterior derivative and (44) computes the
cohomology of the Lie algebra $\widehat{\Theta }\simeq \widetilde{\Theta }.$
The kernel of the first operator in (44) is the space of funtions $%
\widetilde{\mathcal{G}}^{(m-1)}\rightarrow \mathbb{R}$ which are invariant
with respect to the conjugation by the elements of $\widetilde{\mathcal{G}}.$

\begin{proposition}
$\delta :\widehat{\Lambda _{hor}^{m,k}(M^{(m)})}\longrightarrow \widehat{%
\Lambda _{hor}^{m+1,k}(M^{(m+1)})}$
\end{proposition}

Proof: Follows easily from the definitions$.$ \ $\square $

So (44) restricts to $\widehat{}$ -invariant subspaces. We call the
resulting diagram invariant nonlinear horizontal double complex (INHDC) of
the \textit{globalizable} local Lie group $(M,\varepsilon ).$ Proposition 24
implies the following

\begin{corollary}
The restriction of $\delta $ to the kernels of the first horizontal
operators in INHDC computes the Lie group cohomology $H^{\ast }(\widetilde{%
\mathcal{G}},\mathbb{R)}$ of $\widetilde{\mathcal{G}}.$
\end{corollary}

In an attempt to generalize Corollary 25 to higher cohomology groups, we now
observe some further facts about local Lie groups. We recall that if $(M,%
\widetilde{\varepsilon })$ is a local Lie group, that is, if $\widetilde{%
\mathcal{R}}=0,$ then $\widetilde{\Theta }$ is a Lie algebra. So we have the
representation $\mathcal{L}:\widetilde{\Theta }\rightarrow gl(\widetilde{%
\Theta })$ defined by $\mathcal{L}_{\xi }\eta =[\xi ,\eta ],$ where $%
\mathcal{L}$ denotes Lie derivative. More generally, let $\widetilde{%
T_{s}^{r}(M)}$ denote the space of $\widetilde{\varepsilon }$-invariant $%
(r,s)$-tensor fields. We have the representation of $\mathfrak{g=}\widetilde{%
\Theta }$ on $V=$ $\widetilde{T_{s}^{r}(M)}$ defined by

\begin{equation}
\mathcal{L}_{\xi }:\widetilde{T_{s}^{r}(M)}\longrightarrow \widetilde{%
T_{s}^{r}(M)}\text{ \ \ }\xi \in \widetilde{\Theta }
\end{equation}

Observe that for $\xi \in \widehat{\Theta }$, $\mathcal{L}_{\xi }=0$ as an
operator on $\widetilde{T_{s}^{r}(M)}$ since $\widehat{\Theta }$ integrates
to $\widetilde{\mathcal{G}}$ and elements of $\widetilde{T_{s}^{r}(M)}$ are $%
\widetilde{\varepsilon }$-invariant by definition. Now assuming
globalizability, $G=\widehat{\mathcal{G}}$ has a representation $I\mathcal{L}
$ ($I$ denotes integration) on $V=$ $\widetilde{T_{s}^{r}(M)}$ defined by

\begin{equation}
\left( I\mathcal{L}f\right) (\xi )(p)\overset{def}{=}\widehat{\varepsilon }%
(f^{-1}(p),p)_{\ast }\xi (f^{-1}(p))\text{ \ \ \ \ }f\in \widehat{\mathcal{G}%
},\text{ }\xi \in \widetilde{T^{r,s}(M)},\text{ }p\in M
\end{equation}%
and the derivative of the representation (46) is (45), that is, $d(I\mathcal{%
L)=L}$.

Now let $\widehat{H^{r,s}}_{\delta d}$ denote the cohomology group of INHDC
at $(r,s)$ taken first in the horizontal, then in the vertical directions.
Motivated by Corollary 25 and the above general facts, we make (assuming
globalizability) the following conjecture

$\mathbf{C:}$ $\widehat{H^{\ast ,k}}_{\delta d}$ $\simeq H^{k}(G,V)$ where $%
G=\widetilde{\mathcal{G}}$ and $V=\widehat{\Lambda ^{k}(M)}.$

Therefore, if $M$ is compact, $\mathbf{C}$ implies the vanishing of $%
\widehat{H^{\ast ,k}}_{\delta d}$ for $k\geq 1.$ As we indicated above, the
vertical complexes of (40) coincide with the the complex of the composable
cochains in the sense of groupoids with representations as defined in [5].
It is therefore not surprising that for compact $M$ they vanish too by
Proposition 1 in Section 2.1 of [5].

To linearize $\omega \in \Lambda _{x}^{m,k}(M^{m})$, we choose $(m-1)$%
-tangent vectors $\xi ,\eta ,...,\zeta $ at $x$ and let the target variables
$y,z,...,w$ in $\omega (x,y,z,...,w)$ approach $x$ along these directions
\textit{independently }so that

\begin{equation}
L\omega _{I}(x,\xi ,\eta ,...,\zeta )\overset{def}{=}\frac{\partial \omega
_{I}(x,x,...,x)}{\partial y^{a}\partial z^{b}...\partial w^{c}}\xi ^{a}\eta
^{b}...\zeta ^{c}
\end{equation}

Let $\pi :E_{s}\rightarrow M$ be the vector bundle over $M$ whose fiber $\pi
^{-1}(p)$ is the space of $s$-linear maps $T_{p}\times ...\times
T_{p}\rightarrow \mathbb{R}$, that is, $E_{s}=\otimes _{s}T^{\ast }.$ Now
(47) defines a map

\begin{equation}
L:\Lambda _{hor}^{m,k}(M^{m})\longrightarrow \left( \otimes _{m-1}T^{\ast
}\right) \otimes \Lambda ^{k}
\end{equation}

Using (48) we define the linear horizontal double complex (LHDC) and its
invariantization (ILHDC) in such a way that $L$ becomes a homomorphism of
these two biinvariant double complexes. The the $m$' th row of ILHDC
computes $H^{\ast }(\mathfrak{g,}\otimes _{m-1}\mathfrak{g}^{\ast })$ $%
\mathfrak{g=}\widetilde{\Theta }$ and we can show that its vertical
cohomology vanishes for compact $M.$

Finally, consider the covariant differentiation operator $\widehat{\nabla }%
_{X}:T_{s}^{r}(M)\rightarrow T_{s}^{r}(M)$, $X\in \mathfrak{X}(M)=$ the Lie
algebra of smooth vector fields on $M.$ A fundamental fact is expressed by
the formulas

\begin{eqnarray}
\widehat{\nabla }_{\xi } &=&\mathcal{L}_{\eta }\text{ \ \ \ \ \ }\xi \in
\widehat{\Theta },\text{ }\eta =d\Psi (\xi )\in \widetilde{\Theta } \\
\widetilde{\nabla }_{\xi } &=&\mathcal{L}_{\eta }\text{ \ \ \ \ \ }\xi \in
\widetilde{\Theta },\text{ }\eta =d\Psi ^{-1}(\xi )\text{ }\in \widehat{%
\Theta }\text{\ }  \notag
\end{eqnarray}%
where $\Psi $ and $d\Psi $ are defined in Section 2. The formula (49)
continues to be valid if we replace $T_{s}^{r}(M)$ by more general geometric
object bundles as in this note and (49) underlies Propositions 8, 17 and
Proposition 7 in [13]. So in a sense everything in this note and in [13]
reduces to a duality between $\ \widehat{}$ \ and \ $\widetilde{}$ \
together with the concept of invariance on a local Lie group (which, we
believe, is the origin of the concept of torsion), the theory of Lie
derivative on form valued geometric objects (as in [20]), and the relation
between exterior derivative and Lie derivative.

\bigskip

\bigskip

\textbf{References}

\bigskip

\bigskip

[1] E.Abado\u{g}lu, E.Orta\c{c}gil: Intrinsic characteristic classes of a
local Lie group, Portugal. Math. (N.S), Vol. 67 Fasc.4, 2010, 453-483

[2] I.M.Anderson, J.Pohjanpelto: The cohomology of invariant variatonal
bicomplexes, Acta Appl. Math. 41 (1995), 3-19

[3] A.Blaom: Geometric structures as deformed infinitesimal symmetries,
Trans. Amer. Soc. 358 (2006), 3651-3671

[4] C.Chevalley, S.Eilenberg: Cohomology theory of Lie groups and Lie
algebras, Trans. Amer. Math. Soc. 63 (1948), 85-124

[5] M.Crainic: Differentiable and algebroid cohomology, Van Est
isomorphisms, and characteristic classes, Comment. Math. Helv. 78 (2003),
681-721

[6] W.T. Van Est: Group cohomology and Lie algebra cohomology in Lie groups
I, II, Proc Kon. Ned. Akad. 56 (1953), 484-504

[7] G.P.Hochschild, G.D.Mostow: Cohomology of Lie groups, Illinois J. Math.
6 (1962), 367-401

[8] G.Karaali, P.J.Olver, E.Orta\c{c}gil, M.Ta\c{s}k\i n: The symmetry group
of a local Lie group, in progress

[9] I.Kogan, P.J.Olver: The invariant variational bicomplex, Contemp. Math.
285 (2001), 131-144

[10] I.Kolar, P.W.Michor, J.Slovak: Natural Operations in Differential
Geometry, Springer-Verlag, Berlin Heidelberg, 1993

[11] P.J.Olver: Non-associative local Lie groups, J. Lie Theory, 6 (1996),
23-51

[12] E.Orta\c{c}gil: Riemannian geometry as a curved pre-homogeneous
geometry, arXiv.org 1003.3220

[13] E.Orta\c{c}gil: Tensor calculus and deformation theory on a local Lie
group, arxiv.org 1104.5332

[14] T.Tsujishita: On variational bicomplexes associated to differential
equations, Osaka J. Math. 19 (1982), 311-363

[15] A.Verbovetsky: Notes on the horizontal cohomology, Secondary Calculus
and Cohomological Physics (M. Henneaux, I.S.Krasil'shichik, and
A.M.Vinogradov, eds.) Contemporary Mathematics, Amer. Math. Soc.,
Providence, R.I., 1988, 211-231

[16] A.Verbovetsky: Remarks on two approaches to horizontal cohomology:
compatibility complex and the Koszul-Tate resolution, Acta. Appl. Math. 72
(2002), 123-132

[17] A.M.Vinogradov: A spectral sequence associated with a nonlinear
differential equation and algebro-geometric foundations of Lagrangian field
theory with constraints, Soviet Math. Dokl. 19 (1978), 144-148

[18] A.M.Vinogradov: Geometry of nonlinear differential equations, J. Soviet
Math. 17 (1981), 1624-1649

[19] A.Weinstein, P.Xu: Extensions of symplectic groupoids and quantization,
J. Reine Angew. Math. 417 (1991), 159-189

\bigskip \lbrack 20] K.Yano: The Theory of Lie Derivatives and its
Applications, North-Holland, Amsterdam, 1957

\bigskip

\bigskip

Erc\"{u}ment Orta\c{c}gil, Bo\u{g}azi\c{c}i University, Mathematics
Department (emeritus), Bebek, 34342, Istanbul, Turkey

e-mail:ortacgil@boun.edu.tr

\end{document}